\documentclass{amsart}

\usepackage{amssymb}
\usepackage{amsfonts}
\usepackage{amsmath}
\usepackage{amsmath,graphicx,ifthen,fancyhdr}
\usepackage{graphicx}
\newtheorem{theorem}{Theorem}[section]
\newtheorem{lemma}
         {Lemma}[section]
\newtheorem{e-proposition}[theorem]{Proposition}

\newtheorem{e-definition}[theorem]{Definition\rm}

\def\og{\leavevmode\raise.3ex\hbox{$\scriptscriptstyle\langle\!\langle$~}}
\def\fg{\leavevmode\raise.3ex\hbox{~$\!\scriptscriptstyle\,\rangle\!\rangle$}}

\begin{document}
\vspace{-5cm}



%

\vspace{-2.6cm}


\title{Multifractal analysis of multiple ergodic
averages}



\author{Ai-Hua Fan}
\address{LAMFA, UMR 6140 CNRS, Universit\'e de Picardie,
33 rue Saint Leu, 80039 Amiens, France}
\email{ai-hua.fan@u-picardie.fr}

\author{J\"{o}rg Schmeling}
\address{Lund Institute of Technology, Lund University,
Box 118
SE-221 00 Lund, Sweden}
\email{joerg@maths.lth.se}

\author{Meng Wu}
\address{LAMFA, UMR 6140 CNRS, Universit\'e de Picardie,
33 rue Saint Leu, 80039 Amiens, France}
\email{meng.wu@u-picardie.fr}




\begin{abstract}
 In this paper we present a complete solution to the problem of multifractal analysis of multiple ergodic averages in
the case of  symbolic dynamics for functions of two variables depending on the first coordinate.


\vskip 0.5\baselineskip




\end{abstract}

\maketitle




\section{Introduction and results}
Let $T:X\to X$ be a continuous map on a compact metric space $X$. Let $f_1, \cdots, f_\ell$
($\ell\ge 2$) be $\ell$ real valued continuous functions defined on $X$. We consider the following possible limits (for different $x\in X$):
      \begin{equation}\label{MEA}
      M_{f_1, \cdots, f_\ell}(x)=
\lim_{n\rightarrow
\infty}\frac{1}{n}\sum_{k=1}^{n}f_{1}(T^{k}x)f_{2}(T^{2k}x)\cdots
f_{\ell}(T^{\ell k}x).
      \end{equation}
Such limits are widely studied in ergodic theory.
It was proposed in \cite{FLM} to give a multifractal analysis of the multiple ergodic average $M_{f_1, \cdots, f_\ell}(x)$. The authors of \cite{FLM} succeeded in a very special case where $X=\{-1, 1\}^\mathbb{N}$,
$f_k(x) =x_1$ for all $k$ and $T$ is the shift, by using Riesz products. In this note, we shall study the shift map $T$ on the symbolic space $X=\Sigma_m
=S^\mathbb{N}$ with $S=\{0,1, \cdots, m-1\}$ ($m\ge 2$). We assume that $\ell=2$ (the case
$\ell\ge 3$ seems more difficult) and $f_1$ and $f_2$ are H\"{o}lder
continuous. We endow  $\Sigma_m$ with the standard metric: $d(x,y)=m^{-n}$ where $n$ is the largest $k\ge 0$ such that $x_1=y_1$, $\cdots$, $x_k=y_k$.
The Hausdorff dimension of a set $A$ in $\Sigma_m$ will be denoted by $\dim A$.

For any $\alpha \in \mathbb{R}$, define
$$
      L(\alpha)=\{x\in \Sigma_m: M_{f_1, f_2}(x) = \alpha\}.
$$
Let $\alpha_{\min}= \min_{x, y \in \Sigma_m} f_1(x)f_2(y)$ and $\alpha_{\max}=\max_{x, y \in \Sigma_m} f_1(x)f_2(y)$.
Our question is to determine the Hausdorff dimension of $L(\alpha)$.
We further assume that $\alpha_{\min}<\alpha_{\max}$ (otherwise both $f_1$ and
$f_2$ are constant and the problem is trivial).

From classical dynamical system point of view, the set $L(\alpha)$ is not standard and its dimension can not be
described by invariant measures supported on  it. Let us first examine the largest dimension of ergodic measures supported on the set $L(\alpha)$
by introducing the so-called invariant spectrum:
$$
    F_{\rm inv}(\alpha) = \sup \left\{\dim \mu: \mu \ \mbox{\rm ergodic}, \mu(L(\alpha))=1\ \right\}.
$$
Recall that (see \cite{Fan1994})
$$
\dim \mu = \inf\{\dim B: B \ \mbox{Borel \ set},  \mu(B^c)=0 \}.
$$
The dimension $F_{\rm inv}(\alpha)$ is in general smaller than $\dim L(\alpha)$ (compare  the next two theorems). It is even possible that no ergodic
measure is supported on $L(\alpha)$.

\medskip

\begin{theorem}\label{Main1} Let $f_1$ and $f_2$ be two H\"{o}lder continuous functions on $\Sigma_m$.
If $L(\alpha)$ supports an ergodic measure, then
$$
F_{\rm inv}(\alpha)
= \sup \left\{\dim \mu: \mu \ \mbox{\rm ergodic}, \int f_1 d\mu \int f_2 d\mu = \alpha \ \right\}.
$$
    \end{theorem}
    \medskip

    It is known  \cite{FLP} that the above supremum is the dimension of the set of points $x$ such that
    $$
        \lim_{n\to \infty} \frac{1}{n}\sum_{k=1}^n f_1(T^k x)
        \cdot \lim_{n\to \infty} \frac{1}{n}\sum_{k=1}^n f_2(T^k x) =\alpha.
    $$

    Assume that  $f_1$ and $f_2$ are the same function $f$. As a corollary of Theorem~\ref{Main1}, $\mu(L(\alpha))=1$
for some ergodic measure $\mu$
implies $\alpha \ge 0$. So, if $f$ takes a negative  value  $\alpha<0$, then
Theorem 1 shows that there is no ergodic measure supported on $L(\alpha)$. However,
Theorem 2 shows that $\dim L(\alpha)>0$.

 \medskip

In the following we assume that both $f_1$ and $f_2$ depend only on the first coordinate.
For any $s \in \mathbb{R}$, consider the non-linear transfer equation
\begin{equation}\label{transer_equation}
     t_s(x)^2 = \sum_{Ty=x} e^{s f_1(x)f_2(y)}t_s(y).
\end{equation}
It can be proved that the equation admits a unique solution
$t_s: \Sigma_m \to \mathbb{R}_+$,  which depends only on the first coordinate.
Let $dx$ denote the measure of maximal entropy  for the shift on $\Sigma_m$ and let
$$
      P(s) = \log \int_{\Sigma_m} t_s(x) dx + \log m.
$$
Also it can be proved that $P$ is an analytic  convex function and even strictly
convex when $\alpha_{\min} <\alpha_{\max}$ (Lemma \ref{L3}).

\medskip

 \begin{theorem}\label{Main2} Let $f_1$ and $f_2$ be two  functions on $\Sigma_m$ depending only on the first coordinate.
 For  $\alpha \not\in [P'(-\infty), P'(+\infty)]$, we have $L(\alpha)=\emptyset$.
 For  $\alpha \in [P'(-\infty), P'(+\infty)]$, we have
$$\dim
L(\alpha)=\frac{1}{2\log m}\left( P(s_\alpha) - s_\alpha P'(s_\alpha)\right)
$$
where $s_\alpha$ is the unique solution of $P'(s)=\alpha$.
    \end{theorem}
    \medskip

We can prove that  $\alpha_{\min}\le P'(-\infty)\le P'(+\infty)\le \alpha_{\max}$ and that
$\alpha_{\min}= P'(-\infty)$ if and only if there exist $i_0, i_1, \cdots, i_\ell \in S$ ($\ell \ge 1$)
with $i_0=i_\ell$ such that $f_1(i_k) f_2(i_{k+1}) = \alpha_{\min}$ (similar criterion for $\alpha_{\max}= P'(+\infty)$).

Let us look at two examples on $\Sigma_2$.
For $f_1(x) = f_2(x) = 2x_1-1$, we have
$$
     \dim L(\alpha) = \frac{1}{2} + \frac{1}{2} H\left(\frac{1+\alpha}{2}\right),
     \quad F_{\rm inv}(\alpha) = H\left(\frac{1+\sqrt{\alpha}}{2}\right)
$$
where $H(x) = -x\log_2 x -(1-x) \log_2 (1-x)$. See Figure 1  for the
graphs of $\dim L(\alpha)$ and $F_{\rm inv}(\alpha)$.
Remark that $F_{\rm inv}(\alpha) = 0$ but $\dim L(\alpha)>0$ for $-1\le \alpha <0$.
Also remark that $\dim L(\alpha)$ was computed in \cite{FLM} by using Riesz products.
See Figure 2
for the graphs of $\dim L(\alpha)$ and $F_{\rm inv}(\alpha)$ when
$f_1(x) = f_2(x) = x_1$. In the second case $F_{\rm inv}(\alpha)=
H(\sqrt{\alpha})$  and $\dim L(\alpha)$ can be numerically computed
through $P(s)= 2\log t_0(s)$ where $x=t_0(s)$ is the real solution
of the third order algebraic equation
$$
     x^3 -2 x^2 - (e^s -1) x +(e^s -1)=0.
$$

\begin{figure}
\begin{minipage}[t]{0.55\linewidth}
\centering
\includegraphics[width=7.2cm]{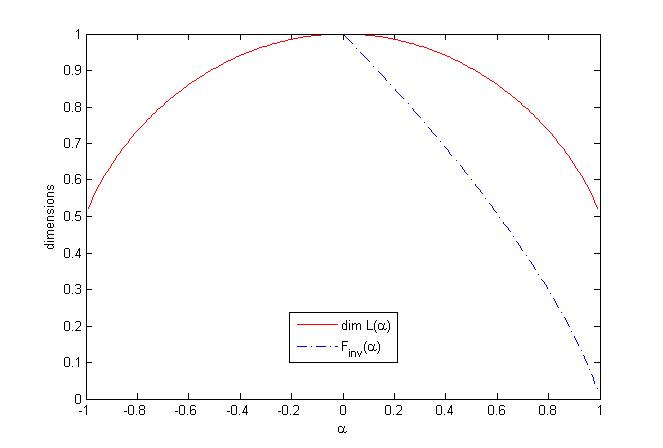}
\caption{\small When $f_1(x) = f_2(x) = 2x_1-1$.} \label{figure1}
\end{minipage}%
\begin{minipage}[t]{0.55\linewidth}
\centering
\includegraphics[width=7.2cm]{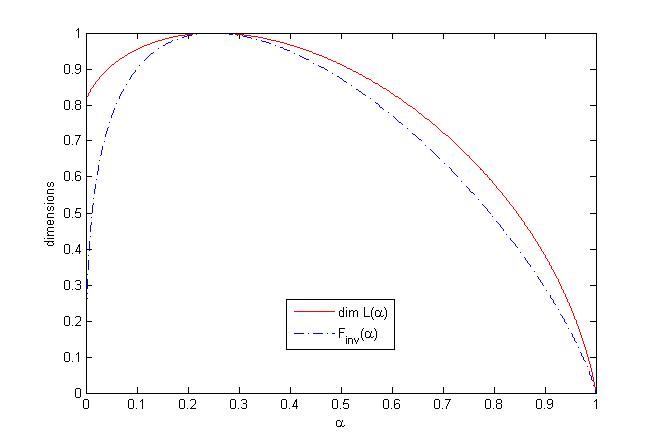}
\caption{\small When $f_1(x) = f_2(x) = x_1$.} \label{figure2}
\end{minipage}
\end{figure}

These two examples show that $F_{\rm inv}(\alpha)< \dim L(\alpha)$
except for some special $\alpha$'s.

The proof of Theorem~\ref{Main2} is based on the following observation.
If $f_1$ and $f_2$ depend only on the first coordinate
$x_1$,  $\sum_k f_1(T^k x)f_2(T^{2k})$ can be decomposed into the sum of $\sum_{j} f_1(T^{i2^j} x)f_2(T^{i2^{j+1}}x)$ with odd $i$, which have independent coordinates.
This observation was used in \cite{FLM} to compute the box dimension
of $X_0=\{x: \forall n, \ x_n x_{2n} = 0  \}$ which is a subset of $L(0)$ (here $f_1(x)=f_2(x) =x_1$ is considered).
The Hausdorff dimension of $X_0$ was later computed in \cite{KPS} where a non-linear transfer operator
characterizes the measure of maximal Hausdorff dimension for $X_0$.

We have stated the results for functions of the form $f_1(x_1) f_2(y_1)$ (product of two functions depending on the first coordinate).
But the results  with obvious modifications hold for functions of the form $f(x_1, y_1)$.

\section{Proof of Theorem~\ref{Main1}}
Let $\mu$ be an ergodic measure such that $\mu(L(\alpha))=1$. Then
$$
   \alpha
   = \lim_{n\to \infty}\frac{1}{n}\sum_{k=1}^n \mathbb{E}_\mu [f_1(T^k x)f_2(T^{2k} x)]
   = \lim_{n\to \infty}\frac{1}{n} \sum_{k=1}^n \mathbb{E}_\mu [f_1( x)f_2(T^{k} x)]
   = \mathbb{E}_\mu [f_1(x) M_{f_2}(x)].
$$
(The first and third equalities are due to Lebesgue convergence theorem and the second one is due to
the invariance of $\mu$). Since $\mu$ is ergodic, $M_{f_2}(x)=\mathbb{E}_\mu f_2$ for $\mu$-a.e. $x$.
So, $\alpha = \mathbb{E}_\mu f_1 \mathbb{E}_\mu f_2$. It follows that
$$
    F_{\rm inv}(\alpha) \le \sup \left\{\dim \mu: \mu \ \mbox{\rm ergodic}, \mathbb{E}_\mu f_1 \mathbb{E}_\mu f_2 = \alpha \ \right\}.
$$
To obtain the inverse inequality, it suffices to observe that the above supremum is attained by a
Gibbs measure $\nu$ which is mixing and that the mixing property implies
$M_{f_1,f_2}(x) = \mathbb{E}_\nu f_1 \mathbb{E}_\nu f_2$ $\nu$-a.e..

\section{Proof of Theorem~\ref{Main2} }


We will prove a result which is a bit more general than Theorem~\ref{Main2}.
Our proof is sketchy and a full proof is contained in \cite{FSW} where other generalizations are also
considered.

Here is the   setting. Let $\varphi: S \times S \to  \mathbb{R}$
be a non constant function with  minimal value $\alpha_{\min}$ and  maximal value $\alpha_{\max}$.
For $\alpha \in \mathbb{R}$, define
$$
     E(\alpha) = \left\{x \in \Sigma_m: \lim_{n \to \infty} n^{-1}\sum_{k=1}^n \varphi(x_k, x_{2k})=\alpha\right\}.
$$

\begin{lemma}\label{L3} For any $s \in \mathbb{R}$, the system
$$
       t_i^2 = \sum_{j=0}^{m-1} e^{s \varphi(i, j)} t_j  \quad (i=0, 1, \cdots, m-1)
$$
admits a unique solution $(t_0(s), t_1(s), \cdots, t_{m-1}(s))$ with strictly positive components, which is an analytic function of $s$.  The function
$$
     P(s)=\log \sum_{j=0}^{m-1} t_j(s)
$$
is strictly convex.
\end{lemma}
\medskip

The proof of the lemma is lengthy. The  existence and uniqueness of the solution are based on the fact that the square roots of right members of the system define an increasing operator on a suitable compact hypercube. The analyticity of the solution is a consequence of the implicit function theorem.

\medskip

\begin{theorem}\label{Main3}
 For any $\alpha \in [P'(-\infty), P'(+\infty)]$, we have
$$
E(\alpha)=\frac{1}{2\log m}\left( P(s_\alpha) - s_\alpha P'(s_\alpha)\right)
$$
where $s_\alpha$ is the unique solution of $P'(s)=\alpha$.
    \end{theorem}

\medskip
The solution  $(t_0(s), t_1(s), \cdots, t_{m-1}(s))$ of the above system allows us to define a Markov measure
$\mu_s$ with initial probability $(\pi(i))_{i\in S}$ and  probability transition matrix
$(p_{i,j})_{S\times S}$  defined by
$$
        \pi(i) = \frac{t_i(s)}{t_0(s) + t_i(s) + \cdots + t_{m-1}(s)},
        \qquad p_{i, j}= e^{s \varphi(i, j)} \frac{t_j(s)}{t_i(s)^2}.
$$

Now decompose the set of positive integers $\mathbb{N}^*$ into
$\Lambda_i$ ($i$ being odd) with $\Lambda_i = \{i 2^k\}_{k\ge 0}$ so that
$\Sigma_m = \prod_{i: 2 \not\mid i} S^{\Lambda_i}$. Take a copy $\mu_s$
on each $S^{\Lambda_i}$ and then define the product measure of these copies. This gives
 a probability measure $\mathbb{P}_s$ on $\Sigma_m$.
Let $\underline{D}(\mathbb{P}_s, x)$ be the lower local dimension
of $\mathbb{P}_s$ at $x$.
\medskip

\begin{lemma} For any $x \in E(\alpha)$, we have
$
       \underline{D}(\mathbb{P}_s, x) \le \frac{1}{2 \log m} [P(s) - \alpha s].
$
\end{lemma}

\medskip
It follows that $\dim E(\alpha) \le \frac{1}{2 \log m} [P(s) - \alpha s]$. Minimizing the right hand side
gives rise to $$
\dim E(\alpha)\le \frac{1}{2 \log m} [P(s_\alpha) - \alpha s_\alpha]$$
where $s_\alpha$ is the solution of $P'(s) =\alpha$. From the lemma, we can deduce that
$L(\alpha)=\emptyset$ if $\alpha \not\in [P'(-\infty), P'(+\infty)]$.
In order to get the inverse inequality, we only have  to show that  $\mathbb{P}_{s_\alpha}$ is supported on $L(\alpha)$.
We first prove the following law of large numbers
by showing the exponential correlation decay of $(F_n)$ under $\mathbb{P}_s$.
\medskip

\begin{lemma} Let $(F_n)$ be a sequence of functions defined on $S\times S$ such that
$\sup_n \sup_{x, y}|F_n(x, y)|<\infty$. For $\mathbb{P}_s$-a.e.  $x \in \Sigma_m$, we have
$$
       \lim_{n\to \infty }\frac{1}{n}\sum_{k=1}^n \left( F_k(x_k, x_{2k})- \mathbb{E}_{\mathbb{P}_s}F_k(x_k, x_{2k}) \right)=0.
$$
\end{lemma}

Applying the above lemma to $F_n(x_n, x_{2n}) = \varphi(x_n, x_{2n})$ for all $n$ and computing $\mathbb{E}_{\mathbb{P}_s}\varphi(x_n, x_{2n})$,
we get
\medskip

\begin{lemma} For $\mathbb{P}_s$-a.e.  $x \in \Sigma_m$, we have
$$
       \lim_{n\to \infty } \frac{1}{n} \sum_{k=1}^n  \varphi(x_k, x_{2k})
       = P'(s).
$$
\end{lemma}

  Thus we  finished the proof for $\alpha \in (P'(-\infty), P'(+\infty))$.
  If $\alpha= P'(-\infty)$ (resp. $P'(+\infty)$),  as in the standard multifractal analysis, we  use
  the probabilities $\mathbb{P}_s$ et let
   $s$ tend to $-\infty$ (resp. $+\infty$).

\vskip1cm

{\em Acknowledgement}  The author would like to thank B. Solomyak for his careful reading of the note
and for the information that the Hausdorff dimension of $L(\alpha)$ is computed in a different way for the special case $m=2, f_1(x)=f_2(x)=x_1$ in \cite{PS}.

\end{document}